\documentclass{llncs}

\usepackage{amsmath,amssymb}
\newcommand{\ie}{i.e.,\ }
\newcommand{\HarmonicSumsP}{\texttt{HarmonicSums}}
\newcommand{\GL}[1]{\textnormal{G}\left(#1\right)}
\renewcommand{\H}[2]{\textnormal{H}_{#1}(#2)}
\newcommand{\HL}[1]{\textnormal{H}_{#1}}
\newcommand{\arxiv}[1]{\texttt{arXiv:#1}}

\allowdisplaybreaks[4]

\begin{document}

\title{Proving two conjectural series for $\zeta(7)$ and discovering more series for $\zeta(7).$}

\author{Jakob Ablinger\thanks{This work was supported by the Austrian Science Fund (FWF) grant SFB F50 (F5009-N15) and has received funding from the European Union’s Horizon 2020 research and innovation programme under the Marie Sk{\l}odowska-Curie grant agreement No. 764850 “SAGEX”.}}

\institute{Research Institute for Symbolic Computation, Johannes Kepler University, Linz, Austria}
\maketitle

\begin{abstract}
We give a proof of two identities involving binomial sums at infinity conjectured by Z-W Sun. In order to prove these identities, we use a recently presented method 
\ie we view the series as specializations of generating series and derive integral representations. 
Using substitutions, we express these integral representations in terms of cyclotomic harmonic polylogarithms. Finally, by applying known relations among the cyclotomic harmonic polylogarithms, we derive the results. 
These methods are implemented in the computer algebra package \HarmonicSumsP.
\end{abstract}

\section{Introduction}\label{sec:Introduction}
In order to prove the two formulas (conjectured in \cite{ZhiWei:2014})
\begin{eqnarray}
\sum_{k=1}^\infty\frac{33 H_k^{(5)}+4/k^5}{k^2\binom{2k}k}&=&-\frac{45}8\zeta(7)+\frac{13}3\zeta(2)\zeta(5)+\frac{85}6\zeta(3)\zeta(4),\label{identity1}\\
\sum_{k=1}^\infty\frac{33 H_k^{(3)}+8/k^3}{k^4\binom{2k}k}&=&-\frac{259}{24}\zeta(7)-\frac{98}9\zeta(2)\zeta(5)+\frac{697}{18}\zeta(3)\zeta(4),\label{identity2}
\end{eqnarray}
where $H_k^{(a)}:=\sum_{i=1}^k\frac{1}{i^a}$, we are going to use a method presented in \cite{PochhammerSums}, therefore we repeat some important definitions and properties (compare~\cite{InvMellin,Ablinger:2014,KauersPaule:2011}).
Let $\mathbb K$ be a field of characteristic~0. A function $f=f(x)$ is called \textit{holonomic} (or \textit{D-finite}) if there exist 
polynomials $p_d(x),p_{d-1}(x),\ldots,p_0(x)\in \mathbb K[x]$  (not all $p_i$
being $0$) such that the following holonomic differential equation holds:
\begin{equation}
 p_d(x)f^{(d)}(x)+\cdots+p_1(x)f'(x)+p_0(x)f(x)=0.
\end{equation}
A sequence $(f_n)_{n\geq0}$ with $f_n\in\mathbb K$ is called
holonomic (or \textit{P-finite}) if there exist polynomials 
$p_d(n),p_{d-1}(n),\ldots,p_0(n)\in \mathbb K[n]$ (not all $p_i$ being $0$) such
that the holonomic recurrence
\begin{equation}
 p_d(n)f_{n+d}+\cdots+p_1(n)f_{n+1}+p_0(n)f_n=0
\end{equation}
holds for all $n\in\mathbb N$ (from a certain point on).
In the following we utilize the fact that holonomic functions are
precisely the generating functions of holonomic sequences: 
for a given holonomic sequence
$(f_n)_{n\geq0}$, the function defined by $f(x) = \sum_{n=0}^{\infty} f_n x^n$ (\ie its 
generating function) is holonomic. 

Note that given a holonomic recurrence for $(f_n)_{n\geq0}$ it is straightforward to 
construct a holonomic differential equation satisfied by its generating function $f(x)=\sum_{n=0}^{\infty} f_n x^n$. For a
recent overview of this holonomic machinery and further literature we
refer to~\cite{KauersPaule:2011}.

In the frame of the proofs we will deal with iterated integrals, hence we define
$$
\GL{f_1(\tau),f_2(\tau),\cdots,f_k(\tau);x}:=\int_0^xf_1(\tau_1)\GL{f_2(\tau),\cdots,f_k(\tau),\tau_1}d\tau_1,
$$
where $f_1(x),f_2(x),\ldots,f_k(x)$ are hyperexponential functions. Note that
$f(x)$ is called \textit{hyperexponential} if $f^\prime(x)/f(x)=q(x),$
where $q(x)$ is a rational function in $x.$

Another important class of iterated integrals that we will come across are the so called cyclotomic harmonic polylogarithms at cyclotomy 3 (compare~\cite{Ablinger:2011te}):
let $m_i:=(a_i,b_i)\in\{(0,0),(1,0),(3,0),(3,1)\}$ for $x\in (0,1)$ we define \textit{cyclotomic polylogarithms at cyclotomy 3}:
\begin{eqnarray}
\H{}{x}&=&1,\nonumber\\
\H{m_1,\ldots,m_k}{x} &=&\left\{ 
		  	\begin{array}{ll}
						\frac{1}{k!}(\log{x})^k,&  \textnormal{if }m_i=(0,0)\\
						  &\\
						\int_0^x{\frac{y^{b_i}}{\Phi_{a_i}(y)} \H{m_2,\ldots,m_k}{y}dy},& \textnormal{otherwise}, 
			\end{array} \right.  \nonumber
\end{eqnarray}
where $\Phi_a(x)$ denotes the $a$th cyclotomic polynomial, for instance $\Phi_1(x) = x - 1 $ and $\Phi_3(x) = x^2 + x + 1.$
We call $k$ the weight of a cyclotomic polylogarithm and in case the limit exists we extend the definition to $x=1$ and write
$$
\HL{m_1,\ldots,m_k}:=\H{m_1,\ldots,m_k}{1}=\lim_{x\to1}\H{m_1,\ldots,m_k}{x}.
$$
Throughout this article we will write $0,1,\lambda$ and $\mu$ for $(0,0),(1,0),(3,0),$ and $(3,1)$, respectively.

Note that cyclotomic polylogarithms evaluated at one posses a multitude of known relations, namely shuffle, stuffle, multiple argument, 
distribution and duality relations, for more details we refer to \cite{Ablinger:2013eba,Ablinger:2011te,MZD}.

\section{Proof of the conjectures}\label{sec:ProofConj}

In order to prove (\ref{identity1}) and (\ref{identity2}) we will apply the method described in \cite{PochhammerSums} and hence we will make use of the 
command \texttt{ComputeGeneratingFunction} which is implemented in the package \HarmonicSumsP\footnote{The package {\tt HarmonicSums} (Version 1.0 19/08/19) together with a Mathematica notebook containing the computations described here can be 
downloaded at \url{http://www.risc.jku.at/research/combinat/software/HarmonicSums}.}\cite{HarmonicSums}. 
Consider the sum left hand side of (\ref{identity1}) and execute (note that in \HarmonicSumsP\ $S[a,k]:=\sum_{i=1}^k\frac{1}{i^a}$)
$$
\textbf{ComputeGeneratingFunction}\left[\frac{33 S[5,k]+4/k^5}{k^2\binom{2k}k},x,\{n,1,\infty \}\right]
$$
which gives (after sending $x\to 1$)
\small
\begin{eqnarray}
 &&\frac{4801781 \text{G}(a,a;1)}{73728}+\frac{451993 \text{G}(0,a,a;1)}{6144}+\frac{10193}{512} \text{G}(0,0,a,a;1)\nonumber\\
 &&+\frac{363}{128} \sqrt{3} \text{G}(a,0,a,a;1)+\frac{1875}{128} \text{G}(0,0,0,a,a;1)+\frac{363}{64}\text{G}(a,a,0,a,a;1)\nonumber\\
 &&+\frac{37}{8} \text{G}(0,0,0,0,a,a;1)+\frac{33}{32} \sqrt{3} \text{G}(a,0,0,0,a,a;1)+\frac{37}{4} \text{G}(0,0,0,0,0,a,a;1)\nonumber\\
 &&+\frac{33}{16} \text{G}(a,a,0,0,0,a,a;1)+\frac{18937121 \text{G}(a;1)}{122880\sqrt{3}}-\frac{895605490019}{5573836800},\label{MathematicaReslut}
\end{eqnarray}
\normalsize
where $0$ represents $1/\tau$ and $a:=\sqrt{\tau}\sqrt{4-\tau}.$\\
Internally \texttt{ComputeGeneratingFunction} splits the left hand side of (\ref{identity1}) into
\begin{eqnarray}\label{Sum1Split}
\sum_{k=1}^\infty x^k \frac{4}{k^7\binom{2k}k}+\sum_{k=1}^\infty x^k \frac{33 H_k^{(5)}}{k^2\binom{2k}k}
\end{eqnarray}
and computes the following two recurrences
\small
\begin{eqnarray*}
 0&=&-(1+k)^7 f(k)+2 (2+k)^6 (3+2 k) f(1+k),\\
 0&=&(1+k)^2 (2+k)^6 f(k)-2 (2+k)^2 (3+2 k) (5+2 k) (55+75 k+40 k^2\\
 &&+10 k^3+k^4) f(1+k)+4 (3+k)^6 (3+2 k) (5+2 k) f(2+k),
\end{eqnarray*}
\normalsize
satisfied by $\frac{4}{k^7\binom{2k}k}$ and $\frac{33 S[5,k]}{k^2\binom{2k}k},$ respectively. 
Then it uses closure properties of holonomic functions to find the following differential equations 
\small
\begin{eqnarray*}
0&=&f(x)+3 (-128+85 x) f'(x)+x (-6906+3025 x) f''(x)\\
&&+14 x^2 (-1541+555 x) f^{(3)}(x)+7 x^3 (-3112+993 x) f^{(4)}(x)\\
&&+42 x^4 (-215+63 x) f^{(5)}(x)+2 x^5 (-841+231 x) f^{(6)}(x)\\
&&+6 x^6 (-23+6 x) f^{(7)}(x)+(-4+x) x^7 f^{(8)}(x),\\
0&=&128 f(x)+8 (-1650+2171 x) f'(x)+2 \left(21870-164445 x+101876 x^2\right) f''(x)\\
&&+2 x \left(264850-761631 x+310438 x^2\right) f^{(3)}(x)\\
&&+4 x^2 \left(354295-599492 x+183087 x^2\right) f^{(4)}(x)\\
&&+2 x^3 \left(694988-826235 x+202454 x^2\right) f^{(5)}(x)\\
&&+8 x^4 \left(76912-70638 x+14483 x^2\right) f^{(6)}(x)\\
&&+x^5 \left(135020-101534 x+17921 x^2\right) f^{(7)}(x)\\
&&+x^6 \left(15020-9614 x+1491 x^2\right) f^{(8)}(x)\\
&&+2 (-4+x) x^7 (-100+31 x) f^{(9)}(x)+(-4+x)^2 x^8 f^{(10)}(x),
\end{eqnarray*}
\normalsize
satisfied by the first and the second sum in (\ref{Sum1Split}), respectively. 

These differential equations are solved using the differential equation solver implemented in~\HarmonicSumsP. 
This solver finds all solutions of holonomic differential equations that can be expressed in
terms of iterated integrals over hyperexponential alphabets~\cite{InvMellin,Ablinger:2014,Bronstein,Singer:99,Petkov:92}; these solutions
are called d'Alembertian solutions~\cite{Abramov:94}, in addition for differential equations of order two it finds all solutions that are Liouvillian~\cite{InvMellinKovacic,Kovacic,Singer:99}.

Solving the differential equations, comparing initial values, summing the two results and sending $x\to 1$ leads to (\ref{MathematicaReslut}).

Since the iterated integrals in (\ref{MathematicaReslut}) only iterate over the integrands $1/\tau$ and $\sqrt{\tau}\sqrt{4-\tau}$ 
we can use the substitution (compare~\cite[Section 3]{BinomialSumIdentities}) $$\tau~\to~(\tau-1)^2/(1+\tau+\tau^2)$$ to compute a representation in terms of cyclotomic harmonic polylogarithms at cyclotomy 3.
This step is implemented in the command \texttt{SpecialGLToH} in \HarmonicSumsP\ and executing this command leads to
\small
\begin{eqnarray*}
&&-3552\HL{\lambda,\lambda,1,1,1,1,1}+1776\HL{\lambda,\lambda,1,1,1,1,\lambda}+3552\HL{\lambda,\lambda,1,1,1,1,\mu}\\
&&+1776\HL{\lambda,\lambda,1,1,1,\lambda,1} -3264\HL{\lambda,\lambda,1,1,1,\lambda ,\lambda}-1776\HL{\lambda,\lambda,1,1,1,\lambda,\mu}\\
&&\hspace{5cm}\vdots\\
&&-1776\HL{\lambda,\lambda,\mu,\mu,\mu,\lambda ,1}+3264\HL{\lambda,\lambda,\mu,\mu,\mu,\lambda,\lambda}+1776\HL{\lambda,\lambda,\mu,\mu,\mu,\lambda,\mu}\\
&&-3552\HL{\lambda,\lambda,\mu,\mu,\mu,\mu,1}+1776\HL{\lambda,\lambda,\mu,\mu,\mu ,\mu,\lambda}+3552\HL{\lambda,\lambda,\mu,\mu,\mu,\mu,\mu},
\end{eqnarray*}
\normalsize
where in total the expression consists of 243 cyclotomic polylogarithms.

Finally, we can use the command \texttt{SpecialGLToH[7,3]} to compute basis representation of the appearing cyclotomic harmonic polylogarithms. \texttt{SpecialGLToH} takes into account shuffle, stuffle, multiple argument, 
distribution and duality relations, for more details we refer to \cite{Ablinger:2013eba,Ablinger:2011te,MZD} and \cite[Section 4]{BinomialSumIdentities}.
Applying these relations we find
\begin{eqnarray}\label{ExampleResult}
-\frac{459}{4} \HL{0,0,1} \HL{\lambda }{}^4-\frac{39}{2} \HL{0,0,0,0,1} \HL{\lambda }{}^2+\frac{45}{8} \HL{0,0,0,0,0,0,1},
\end{eqnarray}
for which it is straightforward to verify that it is equal to the right hand side of~(\ref{identity1}) and hence this finishes the proof. Equivalently we find
\begin{eqnarray*}
\sum_{k=1}^\infty\frac{33 H_k^{(3)}+8/k^3}{k^4\binom{2k}k}&=&\frac{-6273}{20} \HL{\lambda }{}^4 \HL{0,0,1}+49 \HL{\lambda }{}^2 \HL{0,0,0,0,1}+\frac{259}{24} \HL{0,0,0,0,0,0,1},
\end{eqnarray*}
which is equal to the right hand side of (\ref{identity2}).

\section{More identities}\label{sec:Conclusion}
Using the same strategy it is possible to discover also other identities, in the following we list some of the additional identities that we could find:
\begin{eqnarray*}
\sum_{k=1}^\infty\frac{3 H_k^{(2)}-1/k^2}{k^5\binom{2k}k}&=&-\frac{205 \zeta (7)}{18} +\frac{5 \pi ^2 \zeta (5)}{18}+\frac{\pi ^4 \zeta (3)}{18}-\frac{\pi ^7}{486 \sqrt{3}}+\frac{\sqrt{3} c \pi ^3}{8},\\
\sum_{k=1}^\infty\frac{11 H_k^{(3)}+8 H_k^{(2)}/k}{k^4\binom{2k}k}&=&-\frac{7337 \zeta (7)}{216}+\frac{11 \pi ^2 \zeta (5)}{81}+\frac{1417 \pi ^4 \zeta (3)}{4860}-\frac{4 \pi ^7}{729 \sqrt{3}}+\frac{c \pi ^3}{\sqrt{3}},\\
\sum_{k=1}^\infty\frac{2 H_k^{(5)}-H_k^{(3)}/k^2}{k^2\binom{2k}k}&=&-\frac{\zeta (7)}{72}+\frac{8 \pi^2 \zeta (5)}{81}-\frac{17 \pi^4 \zeta (3)}{4860},
\end{eqnarray*}
with $c:=\sum _{i=1}^{\infty } \frac{1}{(3 i+1)^4}+1.$

%

\end{document}